\documentclass[11pt]{article}
\usepackage{amsfonts,amssymb,amsmath}

\usepackage[english]{babel}
\newtheorem{theorem}{Theorem}
\newtheorem{definition}{Definition}

\newtheorem{proof}{Proof}

\newtheorem{proposition}{Proposition}
\newtheorem{lemma}{Lemma}

\newcommand{\beq}{\begin{eqnarray}}
\newcommand{\eeq}{\end{eqnarray}}

\newcommand{\beqt}{\begin{eqnarray*}}
\newcommand{\eeqt}{\end{eqnarray*}}

\newcommand{\be}{\begin{equation}}
\newcommand{\ee}{\end{equation}}

\newcommand{\bl}{\begin{lemma}}
\newcommand{\el}{\end{lemma}}

\newcommand{\bt}{\begin{theorem}}
\newcommand{\et}{\end{theorem}}

\newcommand{\bd}{\begin{definition}}
\newcommand{\ed}{\end{definition}}

\newcommand{\bp}{\begin{proposition}}
\newcommand{\ep}{\end{proposition}}

\newcommand{\bpr}{\begin{proof}}
\newcommand{\epr}{\end{proof}}

\newcommand{\bi}{\begin{itemize}}
\newcommand{\ei}{\end{itemize}}

\newcommand{\ben}{\begin{enumerate}}
\newcommand{\een}{\end{enumerate}}

\newcommand{\Z}{\mathbb Z}
\newcommand{\R}{\mathbb R}
\newcommand{\N}{\mathbb N}

\newcommand{\E}{\mathbb E}

\newcommand{\pee}{\mathbb P}

\newcommand{\eps}{\ensuremath{\epsilon}}

\begin{document}

\title{{\bf A functional limit theorem for a 2d-random walk with dependent marginals}}

\author{Nadine Guillotin-Plantard\footnote{Universit\'e Claude Bernard - Lyon I, institut Camille Jordan,
b\^atiment Braconnier, 43 avenue du 11 novembre 1918,  69622
Villeurbanne Cedex, France. E-mail:
nadine.guillotin@univ-lyon1.fr}, Arnaud Le
Ny\footnote{Universit\'e de Paris-Sud, laboratoire de
math\'ematiques, b\^atiment 425, 91405 Orsay cedex, France.
E-mail: arnaud.leny@math.u-psud.fr}}

\maketitle

\begin{center}
{\bf Abstract}
\end{center}

 We prove a non-standard functional limit theorem
for a two dimensional simple random walk on some randomly oriented
lattices. This random walk, already known to be transient, has
different horizontal and vertical fluctuations leading to
different normalizations in the functional limit theorem, with a
non-Gaussian horizontal behavior. We also prove that the
horizontal and vertical components are not asymptotically
independent.

\medskip
%\footnotesize

\vspace{7cm}



\newpage
\section{Introduction}

The study of random walks on oriented lattices has been recently
intensified with some physical motivations, e.g. in quantum
information theory where the action of a noisy channel on a
quantum state is related to random walks evolving on directed
graphs (see \cite{CP1,CP2}), but they also have their own
mathematical interest. A particular model where the simple random
walk becomes transient on an oriented version of $\mathbb{Z}^2$
has been introduced in \cite{CP2} and extended in \cite{GPLN}
where we have proved a functional limit theorem. In this model,
the simple random walk is considered on an orientation of
$\mathbb{Z}^2$ where the horizontal edges are unidirectional in
some i.i.d. centered random way. This extra randomness yields
larger horizontal fluctuations transforming the usual
normalization in $n^{1/2}$ into a normalization in $n^{3/4}$,
leading  to a non-Gaussian horizontal asymptotic component. The
undirected vertical moves still have standard fluctuations in
$n^{1/2}$ that are thus killed by the larger normalization in the
result proved in \cite{GPLN} (Theorem 4), yielding a null vertical
component in the limit. If these horizontal and vertical
asymptotic components were independent, one could state this
functional limit theorem with an horizontal normalization in
$n^{3/4}$ and a vertical one in $n^{1/2}$, but it might not be the
case. Here, we prove this result  without using independence  and
as a complementary result we indeed prove that these two
asymptotic components are not independent.

%This theorem would be a straightforward extension of Theorem 1 if
%the two components of the limiting process were independent, but
%we also prove here that it is not the case:

%\ben \item RW and statistical mechanics, Quantum, CP1 CP2 GPLN
%LeRoux

%\item horizontally oriented lattices

 %\item CP2, GPLN : transience and FLT

%\item Here: different normalizations, straightforward in case of
%independence, but we also prove the non independence

 %\een

\section{Model and results}

The considered  lattices are oriented versions of $\Z^2$: the
vertical lines are not oriented but the horizontal ones are
unidirectional, the {\em orientation} at a level $y \in \Z$ being
given by a Rademacher random variable $\eps_y= \pm 1$ (say left if
the value is $+1$ and right if
 it is $-1$). We consider here the i.i.d. case where the
 random field $\epsilon=(\epsilon_y)_{y \in \mathbb{Z}}$ has a product
 law $\mathbb{P}_\epsilon=\otimes_{y \in \mathbb{Z}}
 \mathbb{P}_\epsilon^y$ defined on some probability space
 $(A,\mathcal{A},\mathbb{Q})$ with marginals given by $\mathbb{P}^y_\epsilon[\pm 1]= \mathbb{Q}[\epsilon_y = \pm
1]=\frac{1}{2}$.

%We denote by $\mathcal{G}=\sigma\{\epsilon_y, y \in \Z \}$ the
%sub-$\sigma$-algebra of $\mathcal{A}$ generated by the
%orientations.

\bd[Oriented lattices] Let $\eps=(\eps_y)_{y \in \Z}$ be a
sequence of random variables defined as previously. The {\em
oriented lattice} $\mathbb{L}^\eps=(\mathbb{V},\mathbb{A}^\eps)$
is the (random) directed graph with (deterministic) vertex set
$\mathbb{V}=\Z^2$ and (random) edge set $\mathbb{A}^\eps$ defined
by the condition that for $u=(u_1,u_2), v=(v_1,v_2) \in \Z^2$,
$(u,v) \in \mathbb{A}^\eps$ if and only if either $v_1=u_1$ and
$v_2=u_2 \pm 1$, or $v_2=u_2$ and $v_1=u_1+ \eps_{u_2}$. \ed

These orientations will act as generalized random sceneries and
and we denote by $W=(W_t)_{t \geq 0}$ the Brownian motion
associated to it, i.e. such that
\begin{equation}\label{Donsker1}
\Big(\frac{1}{n^{1/2}} \; \sum_{k=0}^{[nt]} \epsilon_{k}\Big)_{t
\geq 0} \;\stackrel{\mathcal{D}}{\Longrightarrow} \big(W_t\big)_{t
\geq 0}.
\end{equation}

%$\lambda$ its characteristic function defined for all $y \in
%\mathbb{Z}$ and  for all $\theta \in \mathbb{R}$ by
%$$
%\lambda(\theta)=\mathbb{Q}\big[e^{i \theta \epsilon_y} \big].
%$$
%We denote by

In this paper,  the notation
$\stackrel{\mathcal{D}}{\Longrightarrow}$ stands for weak
convergence in the space $\mathcal{D}=D([0,\infty[,\R^n)$, for
either $n=1, 2$,  of processes with c\`adl\`ag trajectories
equipped with the Skorohod topology.\footnote{Or sometimes in its
restriction $D([0,T],\R^n)$ for $T >0$. Similarly, we write  $\mathcal{C}$ for the corresponding  spaces of processes with continuous trajectories.}\\

%It is important for us that we have then

%\begin{equation}\label{equiv}
%\lambda(\theta) \sim 1- \frac{\theta^2}{2},\; \; \rm{as} \; \theta
%\to 0.
%\end{equation}

For every realization of $\eps$, one usually means by  simple
random walk on $\mathbb{L}^\eps$ the $\mathbb{Z}^2$-valued Markov
chain $\mathbb{M}=\big(M_n^{(1)},M_n^{(2)} \big)$
  defined on a probability space $(\Omega,
 \mathcal{B},\mathbb{P})$, whose ($\eps$-dependent) transition probabilities are
 defined for all $(u,v) \in \mathbb{V}\times \mathbb{V}$ by
\[\pee[M_{n+1}=v  | M_n=u]=\;\left\{
\begin{array}{lll} \frac{1}{3}  \; &\rm{if} \;  (u,v)   \in \mathbb{A}^\eps&\\
\\
0 \; \; &\rm{otherwise.}&
\end{array}
 \right.
\]
In this paper however, our results are also valid when the
probability of an horizontal move in the direction of the
orientation is $1-p \in [0,1]$ instead of $\frac{1}{3}$, with
probabilities of moving up or down equal thus to $\frac{p}{2}$. We
write  then $m=\frac{1-p}{p}$ for the mean of any geometric random
variable of parameter $p$, whose value is $m=\frac{1}{2}$ in the
standard case  $p=\frac{2}{3}$. We also use a self-similar process $\Delta=(\Delta_t)_{t
\geq 0}$ introduced in \cite{KS} as the asymptotic limit of a
random walk in a random scenery, formally defined for $t \geq 0$
by
$$
\Delta_t= \int_{-\infty}^{+\infty} L_t(x) d W(x)
$$
where $L=(L_t)_{t \geq 0}$ is the local time of a standard
Brownian motion $B=(B_t)_{t \geq 0}$, related to the vertical
component of the walk and independent of $W$. We also denote for
all $t \geq 0$

$$B_t^{(m)}=\frac{1}{\sqrt{1+m}} \cdot B_t\; \; \rm{and} \; \; \Delta_t^{(m)}=\frac{m}{(1+m)^{3/4}} \cdot
\Delta_t.$$

 The following
functional limit theorem has been proved in \cite{GPLN}:

\bt \label{thm1}\cite{GPLN}: \be \label{flt1}
\Big(\frac{1}{n^{3/4}} M_{[nt]} \Big)_{t \geq 0} \;
\stackrel{\mathcal{D}}{\Longrightarrow} \Big(\Delta^{(m)}_t,0 \Big)_{t \geq
0}. \ee \et

We complete here this result with the following theorem: \bt
\label{thm2}:
 \be \label{flt2}
\Big(\frac{1}{n^{3/4}} M^{(1)}_{[nt]}, \frac{1}{n^{1/2}}
M^{(2)}_{[nt]} \Big)_{t \geq 0} \;
\stackrel{\mathcal{D}}{\Longrightarrow} \Big( \Delta_t^{(m)},
B_t^{(m)} \Big)_{t\geq 0}\ee 
and the  asymptotic components $\Delta_t^{(m)}$
and $B_t^{(m)}$ are not independent.

\et

\section{Random walk in generalized random sceneries}

We suppose that there exists some probability space $(\Omega,{\cal
F},\pee)$ on which are defined all the random variables, like e.g.
the orientations $\eps$ and the Markov chain $M$.

\subsection{Embedding of the simple random walk}

We use the orientations  to embed  the 2d-random walk on
$\mathbb{L}^\epsilon$ into two different components: a vertical
simple random walk and an horizontal more sophisticated process.

\subsubsection{Vertical embedding: simple random walk}

The vertical embedding is a one dimensional simple random walk
$Y$, that weakly converges in $\mathcal{D}$ to a standard Brownian
motion $B$:

\be \label{FLT1dY} \Big(\frac{1}{n^{1/2}} Y_{[nt]}\Big)_{t \geq 0}
\; \stackrel{\mathcal{D}}{\Longrightarrow} (B_t)_{t \geq 0}. \ee

The {\em local time} of the walk $Y$ is the discrete-time process
$N(y)=(N_n(y))_{n \in \mathbb{N}}$ canonically defined for all $y
\in \Z$ and $n \in\mathbb{N}$ by

\be \label{ltY}
N_n(y)=\sum_{k=0}^n \mathbf{1}_{Y_k=y} \ee

That is in particular characterized by the following behavior
established in \cite{KS}:

\bl \cite{KS} \label{pre1} $\; \;  \; \; \;\; \;  \; \; \;
\lim_{n\rightarrow\infty}n^{-\frac{3}{4}}\sup_{y\in\Z}N_{n}(y)=0
\; \;  \rm{in} \; \; \mathbb{P}\rm{-probability}.$ \el \noindent

For any reals $a<b$, the fraction of time spent by the process
$\Big(\frac{Y_{[nt]}}{\sqrt{n}}\Big)_{t\geq 0}$ in the
 interval $[a,b)$, during the time interval $\big[0,[nt] \big]$, is defined
 by
 $T_{t}^{(n)}(a,b):=\frac{1}{n}\sum_{a\leq
n^{-\frac{1}{2}}y<b}N_{[nt]}(y)$ 
or
$$\int_{0}^{t} \mbox{\bf 1}_{[a\leq n^{-1/2} Y_{[ns]} < b]}\, ds.$$

One is then particularly interested in analogous quantities for the Brownian motion $(B_{t})_{t\geq 0}$,
i.e.  in a local time $L_{t}(x)$ and in a fraction of time spent in
$[a,b)$ before $t$. If one
defines naturally the former fraction of time to be
$$
\Lambda_{t}(a,b)=\int_0^t \mbox{\bf 1}_{[a \leq  B_{s} < b]}\, ds
$$
then (\cite{LG}) one can define for all $x \in \mathbb{R}$ such a process
$\big(L_t(x)\big)_{t >0}$, jointly continuous in $t$ and $x$, and s.t.,
$$\mathbb{P}-{\rm a.s.}, \;\Lambda_{t}(a,b)=\int_{a}^{b}L_{t}(x)dx.$$

To prove   convergence of the finite-dimensional distributions
in Theorem \ref{thm2}, we need a
more precise relationship between these quantities and  consider
the joint distribution of the fraction of time and the random walk
itself, whose marginals are not necessarily independent.

\bl\label{pre2} For any distinct $t_1, \ldots, t_k \geq 0$ and any $-\infty<a_j<b_j<\infty$ ($j=1,\ldots,k$),

$$
\Big(T_{t_j}^{(n)}(a_j,b_j),\frac{Y_{[n t_j]}}{\sqrt{n}}
\Big)_{1\leq j \leq k} \; \stackrel{\mathcal{L}}{\Longrightarrow}
\; \Big (\Lambda_{t_j}(a_j,b_j),B_{t_j} \Big)_{1\leq j\leq k}
$$

where $\stackrel{\mathcal{L}}{\Longrightarrow}$ means  convergence
in distribution when $n \longrightarrow +
\infty$. \el {\bf Proof:}  For $t\geq 0$, define the projection
$\pi_t$ from ${\cal D}$ to $\R$ as $\pi_t(x) = x_t$. From
\cite{KS}, the map
$$x\in {\cal D} \longrightarrow \int_{0}^{t} \mbox{\bf 1}_{[a\leq x_{s} < b]} \, ds$$
is  continuous on $D([0,T])$ in the Skorohod topology for
any $T\geq t$ for  almost any sample point of the process
$(B_{t})_{t\geq 0}$. Moreover, since almost all paths of the
Brownian motion $(B_{t})_{t\geq 0}$ are continuous at $t$, the map
 $x \rightarrow \pi_t(x)$ is continuous at a.e.
sample points of the process  $(B_{t})_{t\geq 0}$. So, for any
$t\geq 0$, for any $a,b\in\R$ and any $\theta_1\in
\R,\theta_2\in\R$, the map
$$x\in {\cal D} \longrightarrow \theta_1\int_{0}^{t} \mbox{\bf 1}_{[a\leq x_{s} < b]} \, ds + \theta_2 \pi_t(x)$$
 is  continuous  on $D([0,T])$ for any $T\geq t$ at almost all sample points of  $(B_{t})_{t\geq 0}$.
The weak convergence of
$\Big(\frac{Y_{[nt]}}{\sqrt{n}}\Big)_{t\geq 0}$ to the process
$(B_{t})_{t\geq 0}$ implies then the convergence of the law of
$$\sum_{i=1}^{k} \theta_i^{(1)} T_{t_i}^{(n)}(a_i,b_i)+n^{-1/2}\sum_{i=1}^{k} \theta_i^{(2)} Y_{[nt_i]}$$
 to this of $\sum_{i=1}^{k} \theta_i^{(1)} \Lambda_{t_i}(a_i,b_i)+ \sum_{i=1}^{k} \theta_i^{(2)} B_{t_i}$. This proves the lemma using the characteristic function criterion for convergence in distribution.
$\diamond$
\subsubsection{Horizontal embedding: generalized random walk in a
random scenery}

The horizontal embedding is a random walk with geometric jumps:
consider a doubly infinite family $(\xi_i^{(y)})_{i \in
\mathbb{N}^*, y \in \Z}$ of independent geometric random variables
of mean $m=\frac{1-p}{p}$  and define the embedded horizontal
random walk $X=(X_n)_{n\in\N}$ by $X_0=0$ and for $n \geq 1$,
\be \label{GRWRS}
X_n=\sum_{y \in \Z} \eps_y \sum_{i=1}^{N_{n-1}(y)} \xi_i^{(y)}
\ee
with the convention that the last sum is zero when $N_{n-1}(y)=0$.
Define now for $n \in \N$ the random time $T_n$ to be the instant
just after the $n^{{\rm th}}$ vertical move,
\be \label{stoptimes}
T_n=n + \sum_{y \in \Z} \sum_{i=1}^{N_{n-1}(y)} \xi_i^{(y)}.
\ee

Precisely at this time, the simple random walk on
$\mathbb{L}^{\epsilon}$ coincides with its embedding. The
following lemma has been proved in \cite{CP2} and \cite{GPLN}: \bl
\label{lem1} \ben \item $M_{T_n}=(X_n,Y_n),\; \forall n \in \N$.
\item
$$ \frac{T_n}{n} \; \stackrel{n \to \infty}{\longrightarrow}
 1+ m,\;\;\; \mathbb{P}{\rm -almost \; \;surely.}$$

\een \el

\subsection{Random walk in a random scenery}

We call $X$ a generalized random walk in a random scenery
because it is a geometric distorsion of the following {\em random
walk in a random scenery} $Z=(Z_n)_{n \in \N}$ introduced in
\cite{KS} with

$$
Z_n=\sum_{k=0}^n \epsilon_{Y_k} =\sum_{y \in \mathbb{Z}}
\epsilon_y N_{n}(y).
$$

{}From the second expression in terms of the local time of the simple random walk $Y$, it is straightforward to see that its variance is of order $n^{3/2}$, justifying the normalization in $n^{3/4}$ in the functional limit theorem established in \cite{KS}. There, 
the limiting process $\Delta=(\Delta_{t})_{t\geq 0}$ of the
sequence of stochastic processes
$\big({n^{-\frac{3}{4}}Z_{[nt]}}\big)_{t\geq 0}$ is the process
obtained from the random walk in a random scenery when $\Z$ is
changed into $\R$, the random walk $Y$ into a Brownian motion
$B=(B_{t})_{t\geq 0}$ and the random scenery
$(\epsilon_y)_{y\in\Z}$ into a white noise, time derivative in the
distributional sense of a Brownian motion $(W(x))_{x\in\R}$.
Formally replacing $N_{n}(x)$ by $L_{t}(x)$, the  process $\Delta$
can be represented by the stochastic integral
$$\Delta_{t}=\int_{-\infty}^{+ \infty}L_{t}(x)\, dW(x).$$
Since the random scenery is defined on the whole $\Z$ axis, the
Brownian motion $(W(x))_{x\in\R}$ is to be defined with real time.
Therefore, one introduces a pair of independent Brownian motions
$(W_{+}, W_{-})$ so that the limiting process can be rewritten
\begin{equation}\label{4.5}
\Delta_{t}=\int_{0}^{+ \infty}L_{t}(x)\, dW_{+}(x)+\int_{0}^{+
\infty}L_{t}(-x)\, dW_{-}(x).
\end{equation}

In addition to its existence, Kesten and Spitzer have also proved
the

\begin{theorem}\label{KS}\cite{KS}:
$$ \Big(\frac{1}{n^{3/4}} Z_{[nt]} \Big)_{t \geq 0} \;
\stackrel{\mathcal{C}}{\Longrightarrow} (\Delta_t)_{t \geq 0}. $$
\end{theorem}

We complete this result and consider the (non-independent)
coupling between the simple vertical random walk  and the random
walk in  a random scenery and prove:

\begin{theorem}\label{RWIORS}:
$$ \Big(\frac{1}{n^{3/4}} Z_{[nt]},  \frac{1}{n^{1/2}} Y_{[nt]}\Big)_{t \geq 0} \;
\stackrel{\mathcal{D}}{\Longrightarrow} \big(\Delta_t,B_t \big)_{t \geq 0}.
$$
\end{theorem}

%Our main result, Theorem \ref{thm2} will be a consequence of it
%and we mainly focus in the next section on the proof of this
%theorem.

\section{Proofs}
\subsection{Strategy}

The main strategy is to relate the simple random walk on the oriented lattice $\mathbb{L}^\epsilon$ to the random walk in random scenery $Z$ using the embedded process $(X,Y)$. We first prove the  functional limit Theorem \ref{RWIORS} by  carefully carrying the strategy of \cite{KS}, used to prove Theorem \ref{KS}, for a possibly non independent couple $(Z,Y)$. This result extends to the embedded process $(X,Y)$ due to an asymptotic equivalence in probability of $X$ with a multiple of $Z$.  Theorem \ref{thm2} is then deduced from it  using nice convergence properties of the random times (\ref{stoptimes}) and self-similarity. Eventually, we prove that the asymptotic horizontal components of these two-dimensional processes are not independent, using stochastic calculus techniques.

\subsection{Proof of Theorem \ref{RWIORS}}

%We proceed in a standard way, proving a tightness result and
%dealing afterwards with the 

We focus on the convergence of finite dimensional
distributions, because we do not really need the tightness to prove our main result Theorem \ref{thm2}. It could nevertheless be proved in the similar way as the tightness in Lemma \ref{lem6}, see next section.

 \bp\label{pr4.3.1} The finite
dimensional distributions of $\Big(\frac{1}{n^{3/4}} Z_{[nt]},
\frac{1}{n^{1/2}}Y_{[nt]}\Big)_{t\geq 0}$ converge to those of
$(\Delta_{t}, B_{t})_{t\geq 0}$, as $n\rightarrow\infty$. \ep

% We prove the convergence of the finite dimensional distributions
%of $(D_{t}^{n}, n^{-1/2} S_{nt})$ to those of $(\Delta_{t},B_t)$
%using the results of the previous section.

{\bf Proof:} We first identify  the finite dimensional
distributions of $\big(\Delta_{t}, B_{t}\big)_{t\geq 0}$.

\bl\label{pre4} For any distinct $t_{1},\ldots, t_k \geq 0$ and $\theta_{1},\ldots, \theta_k \in\R^2$, the characteristic function
of the corresponding linear combination of $\big(\Delta_{t},B_t \big)$ is given by
$$\E\Big[\exp\Big(i\sum_{j=1}^{k}(\theta_{j}^{(1)}\Delta_{t_{j}}+\theta_{j}^{(2)}B_{t_{j}}) \Big)\Big]=\E\Big[\exp\Big(-\frac{1}{2}\int_{\R}(\sum_{j=1}^{k}\theta_{j}^{(1)}L_{t_{j}}(x))^2\, dx\Big)\exp\Big(i\sum_{j=1}^{k} \theta_{j}^{(2)}B_{t_{j}}  \Big)\Big].$$
\el {\bf Proof:} The function $x\rightarrow
\sum_{j=1}^{k}\theta_{j}^{(1)}L_{t_j}(x)$ being continuous, almost
surely with compact support, for almost all fixed sample of the
random process $(B_t)_{t}$, the stochastic integrals
$$\int_{0}^{+\infty} \sum_{j=1}^{k}\theta_{j}^{(1)}L_{t_j}(x)\, dW_{+}(x)\ \mbox{ and }\ \ \int_{0}^{+\infty} \sum_{j=1}^{k}\theta_{j}^{(1)}L_{t_j}(-x)\, dW_{-}(x)$$
are independent Gaussian random variables, centered, with variance
$$ \int_{0}^{+\infty} \Big(\sum_{j=1}^{k}\theta_{j}^{(1)}L_{t_j}(x)\Big)^2 \, dx\ \mbox{ and }\ \  \int_{0}^{+\infty} \Big(\sum_{j=1}^{k}\theta_{j}^{(1)}L_{t_j}(-x)\Big)^2 \, dx.$$
Therefore, for almost all fixed sample of the random process $B$,
$\sum_{j=1}^{k}\theta_{j}^{(1)}\Delta_{t_{j}}$ is a centered
Gaussian random variable with variance given by
$$\int_{\R} \Big(\sum_{j=1}^{k}\theta_{j}^{(1)}L_{t_j}(x) \Big)^2\, dx.$$
Then we get
\begin{eqnarray*}
& &
\E\Big[\E\Big[e^{i\sum_{j=1}^{k}\theta_{j}^{(1)}\Delta_{t_{j}}}|B_t, t\geq0 \Big]
e^{i\sum_{j=1}^{k}\theta_{j}^{(2)}B_{t_{j}}}\Big]=
\E\Big[e^{-\frac{1}{2}\int_{\R}(\sum_{j=1}^{k}\theta_{j}^{(1)}L_{t_{j}}(x))^2\,
dx } e^{i\sum_{j=1}^{k}\theta_{j}^{(2)}B_{t_{j}}}\Big]. \diamond
\end{eqnarray*}

\smallskip

Hence we have expressed the characteristic function of the linear
combination of $(\Delta_{t},B_t)_{t \geq 0}$ in terms of $B$ and
its local time only. We focus now on the limit of the couple
$\Big(\frac{1}{n^{3/4}} Z_{[nt]}, \frac{1}{n^{1/2}}
Y_{[nt]}\Big)_{t\geq 0}$ when $n$ goes to infinity and introduce
for distinct $t_{j}\geq 0$ and $\theta_{j}\in\R^2$ the
characteristic function

$$\phi_n(\theta_1,\ldots,\theta_{k}):=\E\left[\exp\Big(i n^{-3/4} \sum_{j=1}^{k}
\theta_{j}^{(1)}Z_{[nt_{j}]}\Big)\ \exp\Big(in^{-1/2}\sum_{j=1}^{k}\theta_{j}^{(2)}
 Y_{[nt_j]}\Big)\right].$$
By independence of the random walk $Y$  with  the random scenery
$\epsilon$, one gets
$$\phi_n(\theta_1,\ldots,\theta_{k})=
\E\left[\prod_{x\in\Z}\lambda\Big(n^{-\frac{3}{4}}\sum_{j=1}^{k}\theta_{j}^{(1)}
N_{[nt_{j}]}(x)\Big) \exp\Big(i
n^{-1/2}\sum_{j=1}^{k}\theta_{j}^{(2)}  Y_{[nt_j]}\Big) \right].$$
where $\lambda(\theta)=\mathbb{E}\big[e^{i \theta \epsilon_y} \big]$ is the characteristic function of the orientation $\epsilon_y$, defined for all $y \in \mathbb{Z}$ and for all $\theta \in \mathbb{R}$.
Define now for any $\theta_j\in\R^2$ and $n\geq 1$,
$${\psi}_n(\theta_1,\ldots,\theta_{k}):= \E\left[\exp\Big(-\frac{1}{2}\sum_{x\in\Z}n^{-\frac{3}{2}}
(\sum_{j=1}^{k}\theta_{j}^{(1)}N_{[nt_{j}]}(x))^2\Big) \exp\Big(i
n^{-1/2}\sum_{j=1}^{k}\theta_{j}^{(2)} Y_{[nt_j]}\Big)\right].$$

\bl \label{phipsi} $\; \; \; \; \; \; \; \; \; \;\lim_{n\rightarrow\infty}\Big|\phi_n(\theta_1,\ldots,\theta_{k})-\psi_n(\theta_1,\ldots,\theta_{k})\Big|=0.$
\el

{\bf Proof :} Let $\epsilon >0$ and $A_{n}=\{\omega;
n^{-\frac{3}{4}}\sup_{x\in\Z}|\sum_{j=1}^{k}\theta_{j}^{(1)}N_{[nt_{j}]}(x)|>\epsilon\}$.
Then
\begin{eqnarray*}
& & \left|
\phi_n(\theta_1,\ldots,\theta_{k})-\psi_n(\theta_1,\ldots,\theta_{k})\right|
\;\\ &\leq
&\int_{A_{n}}\left|\prod_{x\in\Z}\lambda\Big(n^{-\frac{3}{4}}\sum_{j=1}^{k}\theta_{j}^{(1)}N_{[nt_{j}]}(x)\Big)
-\exp\Big(-\frac{1}{2}\sum_{x\in\Z}n^{-\frac{3}{2}}(\sum_{j=1}^{k}\theta_{j}^{(1)}N_{[nt_{j}]}(x))^2\Big)\right|d\pee\\
&+&\int_{A_{n}^{c}}\left|\prod_{x\in\Z}\lambda\Big(n^{-\frac{3}{4}}\sum_{j=1}^{k}\theta_{j}^{(1)}N_{[nt_{j}]}(x)\Big)
-\exp\Big(-\frac{1}{2}\sum_{x\in\Z}n^{-\frac{3}{2}}
(\sum_{j=1}^{k}\theta_{j}^{(1)}N_{[nt_{j}]}(x))^2\Big)\right|d\pee.\\
&\leq& 2 \pee(A_{n})+ \int_{A_{n}^{c}}
\left|\prod_{x\in\Z}\frac{\lambda\Big(n^{-\frac{3}{4}}\sum_{j=1}^{k}\theta_{j}^{(1)}N_{[nt_{j}]}(x)\Big)}
{\exp\Big(-\frac{1}{2}\sum_{x\in\Z}n^{-\frac{3}{2}}
(\sum_{j=1}^{k}\theta_{j}^{(1)}N_{[nt_{j}]}(x))^2\Big)}-1\right|d\pee.\\
\end{eqnarray*}
The first term  tends to zero in virtue of
Lemma \ref{pre1}. The second term also vanishes in the limit
because one has $\lambda(\theta)\sim 1-\frac{\theta^2}{2}$ as $|\theta|\rightarrow 0$. Thus Lemma
\ref{phipsi} is proved. $\diamond$\\

The asymptotic behavior of $\phi_n$ will be this of $\psi_n$ and
we identify now its limit with the characteristic function of the
linear combination of $\big( \Delta_t,B_t\big)_{t \geq 0}$ in the
following:

\bl\label{pre3}  For any distinct $t_1,\ldots, t_{k}\geq 0$ and $\theta_{1}, \ldots, \theta_k \in \R^2$, the
distribution of
$$\left(n^{-\frac{3}{2}}\sum_{x\in\Z}\Big(\sum_{j=1}^{k}\theta_{j}^{(1)}N_{[nt_{j}]}(x)\Big)^2,\  n^{-\frac{1}{2}} \sum_{j=1}^{k}
\theta_{j}^{(2)} Y_{[nt_j]}\right)_{j=1 \dots k}$$ converges, as
$n\rightarrow\infty$, to the distribution of
$$\left(\int_{-\infty}^{\infty}\Big(\sum_{j=1}^{k}\theta_{j}^{(1)}L_{t_{j}}(x)\Big)^2\, dx, \sum_{j=1}^{k}
\theta_{j}^{(2)} B_{t_{j}}\right)_{j=1 \dots k}.$$ \el 

{\bf Proof:}  We proceed like in \cite{KS} where a similar result is proved for the horizontal component; although the convergence holds for each component, their possible non-independence prevents to get the convergence for the couple directly and we have to proceed carefully using similar steps and Lemma 2. We decompose the set of all possible
indices into small slices where sharp estimates can be made, and
proceed on them of two different limits on their sizes afterwards.
Define, in a slice of size $a(l,n)=\tau l\sqrt{n}, l\in\Z$, an
average occupation time by
$$
T(l,n)=\sum_{j=1}^{k}\theta_{j}^{(1)}T_{t_{j}}^{(n)}(l\tau,(l+1)\tau)=\frac{1}{n}\sum_{j=1}^{k}\theta_{j}^{(1)}\sum_{a(l,n)\leq y<a(l+1,n)}N_{[nt_{j}]}(y).\\
$$

Define also
$\;\;\;\;\;\;U(\tau,M,n)=n^{-\frac{3}{2}}\sum_{x<-M\tau\sqrt{n}\
\atop\ \mbox{\tiny or}\ x\geq
M\tau\sqrt{n}}(\sum_{j=1}^{k}\theta_{j}^{(1)}N_{[nt_{j}]}(x))^2\;
\; \; \;$ and  $$V(\tau,M,n)=\frac{1}{\tau}\sum_{-M\leq
l<M}(T(l,n))^2+n^{-\frac{1}{2}} \sum_{j=1}^{k}
\theta_{j}^{(2)}Y_{[nt_j]}.$$ Consider
$\delta(l,n)=a(l+1,n)-a(l,n)$ and write
\begin{eqnarray*}
A(\tau, M,n)&:=&n^{-\frac{1}{2}} \sum_{j=1}^{k} \theta_{j}^{(2)} Y_{[nt_j]}+n^{-\frac{3}{2}}\sum_{x\in\Z}\Big(\sum_{j=1}^{k}\theta_{j}^{(1)} N_{[nt_{j}]}(x)\Big)^2-U(\tau,M,n)-V(\tau,M,n)\\
& =&n^{-\frac{3}{2}}\sum_{-M\leq l<M}\sum_{a(l,n)\leq x<a(l+1,n)}
\left(\Big(\sum_{j=1}^{k}\theta_{j}^{(1)}N_{[nt_{j}]}(x)\Big)^2-\frac{n^2\times(T(l,n))^2}{(\delta(l,n))^2}\right).
\end{eqnarray*}

%\medskip

{\it First step:} We first show that $A(\tau, M,n)$ tends in
probability to zero as $n\rightarrow\infty$, for a fixed $\tau$ in
the slice of length $\delta(l,n)$. Fix also $M$ and $n$ and write
\begin{eqnarray*}
& & \E\Big[\Big|
\Big(\sum_{j=1}^{k}\theta_{j}^{(1)}N_{[nt_{j}]}(x)\Big)^2 -
\frac{n^2\times(T(l,n))^2}{(\delta(l,n))^2}\Big|\Big]\\
&=&\E\Big[\Big|\sum_{j=1}^{k}\theta_{j}^{(1)}N_{[nt_{j}]}(x)-\frac{n
\times T(l,n)}{\delta(l,n)}\Big|
\times  \Big|\sum_{j=1}^{k}\theta_{j}^{(1)}N_{[nt_{j}]}(x)+\frac{n \times T(l,n)}{\delta(l,n)}\Big|\Big]\\
&\leq &\E\Big[\Big|\sum_{j=1}^{k}\theta_{j}^{(1)}N_{[nt_{j}]}(x)-\frac{n
\times T(l,n)}{\delta(l,n)}\Big|^2\Big]^{\frac{1}{2}} \times \
\E\Big[\Big|\sum_{j=1}^{k}\theta_{j}^{(1)}N_{[nt_{j}]}(x)+\frac{n
\times T(l,n)}{\delta(l,n)}\Big|^2\Big]^{\frac{1}{2}}.
\end{eqnarray*}
Firstly, $\; \; \; \; \; \E\Big[\Big|\sum_{j=1}^{k}\theta_{j}^{(1)}N_{[nt_{j}]}(x)+ \frac{n
\times T(l,n)}{\delta(l,n)}\Big|^2\Big]$
\begin{eqnarray*}
&\leq& (\delta(l,n))^{-2}\E\Big[\Big(\sum_{j=1}^{k}\sum_{a(l,n)\leq
y<a(l+1,n)}|\theta_{j}^{(1)}|
(N_{[nt_{j}]}(x)+N_{[nt_{j}]}(y))\Big)^2\Big]\\
&\leq &(\delta(l,n))^{-1}\Big(\sum_{j=1}^{k}|\theta_{j}^{(1)}|^2\Big)\  \sum_{j=1}^{k}\sum_{a(l,n)\leq y<a(l+1,n)}\E\Big[(N_{[nt_{j}]}(x)+N_{[nt_{j}]}(y))^2\Big]\\
&\leq&\Big(\sum_{j=1}^{k}|\theta_{j}^{(1)}|^2\Big)\
\sum_{j=1}^{k}\max_{a(l,n)\leq y<a(l+1,n) \atop\ y\neq
x}\E\Big[(N_{[nt_{j}]}(x)+N_{[nt_{j}]}(y))^2\Big]\\
\end{eqnarray*}
\begin{eqnarray*}
&\leq& 2 \Big(\sum_{j=1}^{k}|\theta_{j}^{(1)}|^2\Big)\
\sum_{j=1}^{k}\max_{a(l,n)\leq y<a(l+1,n) \atop\ y\neq
x}\left\{\E[N_{[nt_{j}]}(x)^2]+\E[N_{[nt_{j}]}(y)^2]\right\}\\
&\leq & 2 \Big(\sum_{j=1}^{k}|\theta_{j}^{(1)}|^2\Big)\
\sum_{j=1}^{k}\max_{a(l,n)\leq y<a(l+1,n) \atop\ y\neq
x}\left\{\E[N_{[nt_{j}]}(x)^3]^{2/3}+\E[N_{[nt_{j}]}(y)^3]^{2/3}\right\}
\end{eqnarray*}
and similarly,
 \small
\begin{eqnarray*}
\E\Big[\Big|\sum_{j=1}^{k}\theta_{j}^{(1)}N_{[nt_{j}]}(x)&-&\frac{n
\times T(l,n)}{\delta(l,n)}\Big|^2\Big] \leq
\Big(\sum_{i=1}^{k}|\theta_{j}^{(1)}|^2\Big)\
\sum_{j=1}^{k}\max_{a(l,n)\leq y<a(l+1,n) \atop\ y\neq
x}\E\Big[\big(N_{[nt_{j}]}(x)-N_{[nt_{j}]}(y)\big)^2\Big].
\end{eqnarray*}
\normalsize
 Thus, using Lemma 1 and 3 from \cite{KS}, we have for large $n$,
$$\E\Big[\Big|A(\tau,M,n)\Big|\Big] \leq  C (2M+1)\tau^{3/2}.$$

We will afterwards consider the limit $M\tau^{3/2}$ goes to zero to
approximate the stochastic integral of the local time $L_t$, and
this term will then go to zero. Moreover, we have
\begin{eqnarray*}
\pee[U(\tau,M,n)\neq 0]
&\leq &\pee[N_{[nt_{j}]}(x)>0\ \mbox{for some}\ x\ \mbox{such that}\ |x|>M\tau\sqrt{n}\ \mbox{and}\ 1\leq j\leq k]\\
&\leq &\pee\Big[N_{\max([nt_{j}])}(x)>0\ \mbox{for some}\ x\
\mbox{such that}\
|x|>\frac{M\tau}{\sqrt{\max(t_{j})}}\sqrt{\max([nt_{j}])}\; \Big].
\end{eqnarray*}

{}From item b) of Lemma 1 in \cite{KS} , we can choose $M\tau$ so
large that $\pee\big[U(\tau,M,n)\neq 0\big]$ is small. Then, we
have proved that for each $\eta>0$, we can choose $\tau,M$ and
large $n$ such that
$$\pee\Big[\Big|n^{-\frac{1}{2}} \sum_{j=1}^{k} \theta_{j}^{(2)} Y_{[nt_j]}+n^{-\frac{3}{2}}\sum_{x\in\Z}(\sum_{j=1}^{k}
\theta_{j}^{(1)}N_{[nt_{j}]}(x))^2-V(\tau,M,n)\Big|>\eta\Big]\leq
2\eta.$$ \noindent
\medskip
{\it Second step:} From Lemma \ref{pre2}, $V(\tau,M,n)$ converges
in distribution, when $n\rightarrow\infty$, to
$$\frac{1}{\tau}\sum_{-M\leq  l<M}\Big(\sum_{j=1}^{k}\theta_{j}^{(1)}\int_{l\tau}^{(l+1)\tau}L_{t_{j}}(x)dx\Big)^2+\sum_{j=1}^{k}\theta_j^{(2)} B_{t_j}.$$
The function $x\rightarrow L_{t}(x)$ being continuous and having a.s
compact support,
$$\frac{1}{\tau}\sum_{-M \leq l<M}\Big(\sum_{j=1}^{k}\theta_{j}^{(1)}\int_{l\tau}^{(l+1)\tau}L_{t_{j}}(x)dx\Big)^2+\sum_{j=1}^{k}\theta_j^{(2)} B_{t_j}$$
converges, as $\tau\rightarrow 0, M\tau\rightarrow\infty$, to
$$\int_{-\infty}^{\infty}(\sum_{j=1}^{k}\theta_{j}^{(1)}L_{t_{j}}(x))^2\,   dx+\sum_{j=1}^{k}\theta_j^{(2)} B_{t_j}.\; \diamond$$

 Putting together Lemma \ref{pre4}, \ref{phipsi}  and
\ref{pre3} gives Proposition \ref{pr4.3.1}, that proves Theorem
\ref{RWIORS}. $\diamond$

% From Lemma \ref{pre3},
%\begin{eqnarray*}
%& &\lim_{n\rightarrow\infty}\psi_n(\theta_1,\ldots,\theta_{k})\\
%&=&\E\left[\exp\Big(-\frac{1}{2}\int_{-\infty}^{\infty}(\sum_{j=1}^{k}\theta_{j}^{(1)}
%L_{t_{j}}(x))^2\, dx\Big) \exp\Big(i\sum_{j=1}^{k} \theta_{j}^{(2)} B_{t_{j}}\Big)\right]\\
%&=&\E\left[\exp\Big(i\sum_{j=1}^{k}\theta_{j}^{(1)}\Delta_{t_{j}}\Big)\
%\exp\Big(i\sum_{j=1}^{k} \theta_{j}^{(2)} B_{t_{j}}\Big)\right],
%\mbox{by Lemma \ref{pre4}}.
%\end{eqnarray*}
%Then, the proposition is proved.

\subsection{Proof of Theorem \ref{thm2}}

We get the convergence of Theorem \ref{thm2}  from Theorem \ref{RWIORS} and Lemma \ref{lem1} and focus first on the embedded process $(X,Y)$:

\bl \label{lem6}
$$ \Big(\frac{1}{n^{3/4}} X_{[nt]},  \frac{1}{n^{1/2}} Y_{[nt]}\Big)_{t \geq 0} \;
\stackrel{\mathcal{D}}{\Longrightarrow} \Big(m \cdot
\Delta_t,B_t\Big)_{t \geq 0}.
$$
\el

{\bf Proof}: We first prove the tightness of the family. The
second component is  tight in ${\cal D}$ (see Donsker's theorem in
\cite{Bil}), so to prove the proposition we only have to prove the
tightness of the first one in ${\cal D}$.
By Theorem 13.5 of Billingsley \cite{Bil}, it suffices to prove
that there exists $K>0$ such that for all $t,t_1,t_2\in[0,T],
T<\infty,$ s.t. $t_{1}\leq t\leq t_{2},$ for all $n\geq 1$,
\begin{equation}\label{pro}
\E\Big[|X_{[nt]}-X_{[nt_1]}| \cdot \ |X_{[nt_2]}-X_{[nt]}|\Big]\leq K
n^{3/2} |t_{2}-t_{1}|^{\frac{3}{2}}.
\end{equation}
Using Cauchy-Schwarz inequality, it is enough to prove that there
exists $K>0$ such that for all $t_1\leq t$, for all $n\geq 1$,
\begin{equation}\label{pro1}
\E\Big[|X_{[nt]}-X_{[nt_1]}|^2\Big]\leq K
n^{3/2} |t_{2}-t_{1}|^{\frac{3}{2}}.
\end{equation}
Since the $\eps'$s are independent and centered, we have
$$
\E\Big[|X_{[nt]}-X_{[nt_1]}|^2\Big]=\sum_{x\in\Z}
\E\Big[\sum_{i=N_{[nt_1]-1}(x)+1}^{N_{[nt]-1}(x)}\sum_{j=N_{[nt_1]-1}(x)+1}^{N_{[nt]-1}(x)}
\E[\xi_{i}^{(x)}\xi_{j}^{(x)}| Y_k, k\geq 0]\Big].
$$
{}From the inequality
$$0\leq \E[\xi_{i}^{(x)}\xi_{j}^{(x)}]\leq m^2 + \mbox{Var} (\xi_{i}^{(x)})=C,$$
we deduce that
\begin{eqnarray*}
\E\Big[|X_{[nt]}-X_{[nt_1]}|^2 \Big]&\leq & C
\sum_{x\in\Z}\E\Big[(N_{[nt]-1}(x)-N_{[nt_{1}]-1}(x))^2\Big]= C
\sum_{x\in\Z}\E\Big[(N_{[nt]-[nt_{1}]-1}(x))^2\Big].
\end{eqnarray*}

{}From item d) of Lemma 1 in \cite{KS}, as $n$ tends to infinity,
$$ \E\big[\sum_x N_n^{2} (x)\big] \sim C n^{3/2},$$
and there exists some constant $K>0$ such that
\begin{eqnarray*}
\E\Big[|X_{[nt]}-X_{[nt_1]}|^2\Big]\leq  K
\Big([nt]-[nt_{1}]-1\Big)^{\frac{3}{2}}\leq K
n^{\frac{3}{2}}\Big(t-t_1\Big)^{\frac{3}{2}}.
\end{eqnarray*}
We get the tightness of the first component by dividing $X_n$  by $n^{3/4}$, and eventually the tightness of the properly normalized embedded process.\\

To deal with finite dimensional distributions, we rewrite $X_n=X_n^{(1)}+m Z_{n-1}$ with 
$$
X_n^{(1)}=\sum_{y \in \Z} \epsilon_y
\sum_{i=1}^{N_{n-1}(y)} \big( \xi_i^{(y)} -m \big).
$$
Using the $\mathbb{L}^2-$convergence proved in the proof of Proposition 2 in
\cite{GPLN},
$$
\frac{X_n^{(1)}}{n^{3/4}} \stackrel{n \to \infty}{\longrightarrow} 0, \, \,   {\rm in\  Probability}
$$
one gets that the finite dimensional distributions of
$\Big(\frac{X_{[nt]}}{n^{3/4}}, \frac{Y_{[nt]}}{n^{1/2}} \Big)_{t\geq 0}$ are asymptotically
equivalent to those of $\Big(m \cdot \frac{Z_{[nt]}}{n^{3/4}}, \frac{Y_{[nt]}}{n^{1/2}} \Big)_{t\geq 0}$. One concludes then using Theorem 4. $\diamond$\\

In the second step of the proof of Theorem \ref{thm2}, we use Lemma 3 of
\cite{GPLN} and  that
$M_{T_n}=\big(M_{T_n}^{(1)},M_{T_n}^{(2)}\big)=(X_n, Y_n)$ for any
$n$ with

$$
\frac{T_n}{n} \; \stackrel{n \to \infty}{\longrightarrow} 1+m \;,
\; \;\mathbb{P}-{\rm almost \; \; surely}
$$
and the self-similarity of the limit process $\Delta$ (index $3/4$) and of the Brownian motion $B$ (index $1/2$). Using the strict increasing of $(T_n)_{n \in
\N}$, there exists a sequence of integers $(U_n)_n$ which tends to
infinity and such that $ T_{U_n} \leq n < T_{U_{n}+1}$. More
formally, for any $n\geq 0$, $U_n=\sup\{k\geq 0; T_k\leq n\}$,
$(U_{[nt]}/n)_{n\geq 1}$ converges a.s. to the continuous function
$\phi(t):=t/(1+m)$, so from Theorem 14.4 from \cite{Bil},
$$\Big(\frac{1}{n^{3/4}} M^{(1)}_{T_{U_{[nt]}}}, \frac{1}{n^{1/2}}
M^{(2)}_{T_{U_{[nt]}}} \Big)_{t \geq 0} \;
\stackrel{\mathcal{D}}{\Longrightarrow}  \big( m\Delta_{\phi(t)},
B_{\phi(t)}\big)_{t \geq 0}.$$
Using Lemma 7, the processes  $(m\Delta_{\phi(t)}, B_{\phi(t)} )_t$ and
$\Big(\frac{m}{(1+m)^{3/4}} \Delta_t, \frac{1}{(1+m)^{1/2}} B_t \big)_t$ have the same law, so
$$\Big(\frac{1}{n^{3/4}} M^{(1)}_{T_{U_{[nt]}}}, \frac{1}{n^{1/2}}
M^{(2)}_{T_{U_{[nt]}}} \Big)_{t \geq 0} \;
\stackrel{\mathcal{D}}{\Longrightarrow} \big( \Delta_t^{(m)},
B_t^{(m)}\big)_{t \geq 0}$$ with
$\Delta_t^{(m)}=\frac{m}{(1+m)^{3/4}} \; \cdot \Delta_t$ and
$B_t^{(m)}=\frac{1}{\sqrt{1+m}} \cdot B_t$ for all $t \geq 0$.
Now, $M_{[nt]}^{(2)}=M^{(2)}_{T_{U_{[nt]}}}$ and $M^{(1)}_{[nt]} =
M^{(1)}_{T_{U_{[nt]}}} +\big( M^{(1)}_{[nt]} -
M^{(1)}_{T_{U_{[nt]}}}\big)$, so
$$ \Big|M^{(1)}_{[nt]} - M^{(1)}_{T_{U_{[nt]}}}\Big| \leq \Big|M^{(1)}_{T_{U_{[nt]}+1}} - M^{(1)}_{T_{U_{[nt]}}}\Big|
= \xi_{N_{U_{[nt]}}(Y_{U_{[nt]}})}^{(Y_{U_{[nt]}})}.$$
By remarking that for every $T>0$,
$$ \pee\Big[ \sup_{t\in[0,T]} \frac{1}{n^{3/4}} \xi_{N_{U_{[nt]}}(Y_{U_{[nt]}})}^{(Y_{U_{[nt]}})}\geq \eps\Big]\leq [nT] \cdot \pee [ \xi_1^{(1)} \geq \eps n^{3/4}]  \leq \frac{[nT] \E[|\xi_1^{(1)}|^2]}{\eps^2 n^{3/2}}=o(1),$$
we deduce that for any $T>0$, $\Big(\frac{M^{(1)}_{[nt]} -
M^{(1)}_{T_{U_{[nt]}}}}{n^{3/4}}, 0\Big)_{t\in [0,T]}$ converges as
an element of ${\cal D}$ in $\mathbb{P}$-probability to 0.
Finally, we get the result:
$$
\Big(\frac{1}{n^{3/4}} M^{(1)}_{[nt]}, \frac{1}{n^{1/2}}
M^{(2)}_{[nt]} \Big)_{t \geq 0} \;
\stackrel{\mathcal{D}}{\Longrightarrow} \big(
\Delta_t^{(m)},B_t^{(m)}\big)_{t \geq 0}.$$

Let us prove now that we could not deduce this result from the convergence of the components because the limiting horizontal and vertical  components
 are not  independent. It is enough to prove that $\Delta_1$ and $B_1$ are not
independent and we use that conditionnally to $(B_t)_{0 \leq t \leq 1}$,
 the random variable $\Delta_1$ is the sum of the stochastic integrals
$$\int_{0}^{+\infty} L_{1}(x)\, dW_{+}(x)\ \mbox{ and }\ \
\int_{0}^{+\infty} L_{1}(-x)\, dW_{-}(x)$$
which are independent Gaussian random variables, centered, with variance
$$ \int_{0}^{+\infty} L_{1}(x)^2 \, dx\ \mbox{ and }\
\  \int_{0}^{+\infty} L_{1}(-x)^2 \, dx.$$

Denote by $ V_1:=  \int_{\R} L_1^2(x) dx$ the self-intersection
time of the Brownian motion $(B_t)_{t\geq 0}$ during the time
interval $[0,1]$.

\bl For $n \in \mathbb{N}$ even, one has $
\mathbb{E}\big[V_1 \cdot B_1^n  \big]=C(n) \cdot
\mathbb{E}[B_1^n]$: $V_1$ and $B_1$ are not
independent. \el

 {\bf Proof:} For every $x\in\R$, define $J_{\varepsilon} (x)=\frac{1}{2\varepsilon}
\int_{0}^{1} {\mbox 1}_{\{|B_s -x|\leq \varepsilon\}}\, ds$. Then,
$L_1^2(x)$ is the almost sure limit of $\big(J_{\varepsilon}
(x)\big)^2$ as $\varepsilon\rightarrow 0$ so that
$$
V_1 \cdot B_1^n = \int_\R \Big(\lim_{\varepsilon\rightarrow 0}
J_{\varepsilon}(x)^2 B_1^n \Big) dx
$$
and by Fubini's theorem for $n \in \N$ even,
$$
\E\big[ V_1 \cdot B_1^n \big] = \int_{\R}\E\left[ \lim_{\varepsilon\rightarrow 0}
J_{\varepsilon}(x)^2 B_1^n \, \right] dx.
$$
{}From the occupation times formula, for every $x\in \R$, for every $\varepsilon >0$,
$$J_{\varepsilon} (x)\leq L_1^{*}:= \sup_{x\in \mathbb{R}} L_1(x).$$
So, for every $x\in \R$, for every $\varepsilon >0$, $J_{\varepsilon}(x)^2 B_1^n$ is dominated by $(L_1^{*})^2\ B_1^n$ which belongs to $\mathbb{L}^1$ since $L_1^{*}$ and  $B_1$ have moments of any order (see \cite{Rev} for instance). By dominated convergence theorem, we get
$$
\mathbb{E}\big[ V_1 \cdot B_1^n \big] =\int_{\R}
\lim_{\varepsilon\rightarrow 0} \E \left[
J_{\varepsilon}(x)^2 B_1^n \, \right] dx.
$$
 But, when $(p_t)_t$ is the
Markov transition kernel of the Brownian motion $B$,
\begin{eqnarray*}
\E \left[
J_{\varepsilon}(x)^2 B_1^n \, \right]& =&\frac{1}{2\varepsilon^2}\E\left[ \int_{0<s<t\leq 1}
{\mbox 1}_{\{|B_s -x|\leq \varepsilon\}}{\mbox 1}_{\{|B_t -x|\leq
\varepsilon\}}\, B_1^n dsdt\right]\\
&=&\frac{1}{2\varepsilon^2} \int_{\R^3}
\int_{0<s<t\leq 1} {\mbox 1}_{\{|y -x|\leq \varepsilon\}}{\mbox 1}_{\{|z -x|\leq \varepsilon\}} p_s(0,y) p_{t-s}(y,z)
p_{1-t}(z,u) u^n \, ds dt dy dz du \\
&=& 2 \int_{\R} du \int_{0<s<t\leq 1} ds dt
\left[\frac{1}{4\varepsilon^2}
\int_{x-\varepsilon}^{x+\varepsilon}
\int_{x-\varepsilon}^{x+\varepsilon} p_s(0,y) p_{t-s}(y,z)  p_{1-t}(z,u) u^n \, dy dz \right]\\
\end{eqnarray*}
which converges as $\varepsilon\rightarrow
0$ to
$$\int_{\R} p_s(0,x) p_{t-s}(x,x)p_{1-t}(x,u) u^n\, du.$$
 We deduce that
\begin{eqnarray*}
\E[V_1 \cdot B_1^n]&=& 2 \int_{0<s<t\leq 1} p_{t-s}(0,0) ds dt\int_{\R}\left[ \int_{\R} p_s(0,x) p_{1-t}(x,u)\, dx\right] u^n du\\
&=& 2 \int_{0<s<t\leq 1} p_{t-s}(0,0) \left[\int_{\R} p_{1-t+s}(0,u)u^n \, du\right] ds dt.\\
\end{eqnarray*}
Now, by the scaling property of the Brownian motion,
$$\int_{\R} p_{1-t+s}(0,u)u^n \, du= \E[ B_{1-t+s}^n] = (1-t+s)^{n/2} \E[B_1^n].$$
Therefore, $\E[V_1\cdot B_1^n]= C(n)\cdot \E[B_1^n]$
where
\begin{eqnarray*}
C(n) &=& 2 \int_{0<s<t\leq 1} \frac{(1-t+s)^{n/2}}{\sqrt{2\pi
(t-s)}} \, ds dt. \; \diamond
\end{eqnarray*}

To get the non-independence, one  computes then for $n$ even
\begin{eqnarray*}
\mathbb{E}\big[ \Delta_1^2 \cdot B_1^n \big]
&=& \mathbb{E}\big[\mathbb{E}\big[\Delta_1^2 | B_s, 0\leq s\leq 1\big]\cdot B_1^n \big] = \mathbb{E}\big[V_1^2 \cdot B_1^n \big]\\
&\neq& \mathbb{E}\big[V_1^2 \big]
\cdot\mathbb{E}\big[ B_1^n \big]
=\mathbb{E}\big[\Delta_1^2 \big]\cdot\mathbb{E}\big[ B_1^n \big]
\end{eqnarray*}
leading to the non-independence of $\Delta_1$ and $B_1$.
$\diamond$.

\section{Conclusions and open questions}

The functional limit theorem we have proved here, with an
horizontal component normalized in $n^{3/4}$ with a non-Gaussian
behavior and a more standard vertical one normalized in $n^{1/2}$,
strongly indicates the possibility of a local limit theorem where
the probability for the walk to come back at the origin would be
of order $n^{5/4}$, in complete coherence with the transience
result of \cite{CP2}. This result is not straightforward and an extra
work is needed to get it; this is a work in progress.
Other interesting questions could concern different lattices, with
e.g. also vertical orientations, but the peculiarities of these
studies on oriented lattices is that the methods used are not
robust for the moment. Getting robust methods for general oriented
lattices seems to be a nice challenge.\\

{\bf Acknowledgments :}  We thank D. P\'etritis and W. Werner for having independently raised
this question, S. Friedli, Y. Le Jan and O. Raimond for their
interest in the independence question, and A. Lachal and J.F. Le
Gall for their disponibility to answer to local times questions.

\end{document}